\documentclass[preprint]{elsarticle}

\usepackage{lineno,hyperref}

\usepackage{graphicx}
\usepackage{amsmath,amssymb}
\usepackage{verbatim}
\usepackage{mathptmx}      
\usepackage{latexsym}
\usepackage{marvosym}
\modulolinenumbers[5]

\journal{Journal of \LaTeX\ Templates}

\newtheorem{thm}{Theorem}

\newdefinition{rmk}{Remark}
\newproof{pf}{Proof}
\newproof{pot}{Proof of Theorem \ref{thm2}}
\newcommand{\SP}{\mathop{\mathrm{supp}}\nolimits}

\begin{document}

\begin{frontmatter}

\title{Constructive description of H{\"o}lder-like classes on an arc in $\mathbb{R}^3$ by means of harmonic functions}

\author[rvt]{Tatyana~A.~Alexeeva\corref{cor1}}
\ead{tatyanalexeeva@gmail.com}

\author[rvt,focal]{Nikolay~A.~Shirokov}
\ead{nikolai.shirokov@gmail.com}

\cortext[cor1]{Corresponding author}
\address[rvt]{Department of Mathematics, 

National Research University Higher School of Economics, 

3A Kantemirovskaya Ul., St.~Petersburg, 194100, Russia}


\address[focal]{Department of Mathematical Analysis, Faculty of Mathematics and Mechanics, 

St.~Petersburg State University, 28 Universitetsky prospekt, Peterhof, St.~Petersburg, 198504, Russia}




\fntext[myfootnote]{The second author was supported by the RFBR grnt 17-01-00607}

\begin{abstract}
We give a constructive description of H{\"o}lder-like classes of functions on  chord-arc curves in $\mathbb{R}^3$ in terms  of a rate of approximation by harmonic functions in shrinking neighborhoods of those curve.
\end{abstract}


\end{frontmatter}

\linenumbers

\section{Introduction}

The   constructive description of classes of functions in terms of a rate of approximation   by functions taken from  specific sets  (polynomials, rational functions, entire functions, etc.) was initiated by  D.~Jackson and S.~N.~Bernstein at the beginning of the 20th century. Nevertheless, a natural   problem of a constructive description of H{\"o}lder classes on a segment in terms of a rate of approximation by algebraic polynomials was solved only in 1956 \cite{Dzyadyk56}, [\cite{Dzyadyk08},~ch.~7].
Since then problems concerning  constructive description of classes of functions defined on domains in the complex plane have played a central role in   approximation theory. Many authors were involved in the following problem: let $G$ be a Jordan region in the complex plane $\mathbb{C}$, and let $H(G)$ be a class of functions $f$ analytic in the interior $\stackrel{\circ}{G}$ of $G$ and continuous (or smooth)  on the closure of $G$. What is the scale of approximation of functions from $H(G)$ by algebraic polynomials which makes it possible to find  the rate of smoothness of relevant functions? V.~K.~Dzyadyk (\cite{Dzyadyk59,Dzyadyk62,Dzyadyk63,Dzyadyk_Nikol63}) introduced a special type of weights $\rho_{1/n}(z)$ on the boundary $\Gamma$ of   $G$ such that the condition that $f$ is   analytic  in  $\stackrel{\circ}{G}$ and satisfies the $\alpha$-H{\"o}lder condition, $\alpha > 0$, $\alpha\notin \mathbb{N}$, is equivalent to the possibility of approximating $f$ by polynomials $P(z)$ of degree $\le n$ with the property 

\begin{equation*}
|f(z)-P_n(z)|\le C_f \rho_{1/n}^{\alpha}(z),\, z\in\Gamma . \eqno(\star) 
\end{equation*}

So,  for various regions in $\mathbb{C}$, the weights $\rho_{1/n}^{\alpha}(z)$ were a successful scale for a constructive description of the above-mentioned classes of functions. The main problem in that direction was to weaken the assumptions concerning the boundary $\Gamma$. The results progressed from a piecewise smooth in some sense \cite{Dzyadyk66,Lebedev71,Shirokov72} to a chord-arc \cite{Dzyadyk77} and finally to a quasiconformal property of a Jordan curve $\Gamma$ \cite{Belyi77}.

If turned out that if a function $f$ can be approximated by polynomials $P_n(z)$ of degree $\le n$ as in ($\star$), then $f$ is analytic in $\stackrel{\circ}{G}$ and satisfies the $\alpha$-H{\"o}lder condition for any Jordan domain $G$ \cite{Lebedev66,LebedevT70}.

In the case where the boundary $\Gamma = \partial G$ has  cusps,  the polynomial approximation with the rate $const\cdot \rho_{1/n}^{\alpha}(z)$ is appropriate not for all functions satisfying the $\alpha$-H{\"o}lder condition \cite{Shirokov01,Andrievskii85}. Consequently,  in the case of an arbitrary Jordan region, the scale $\rho_{1/n}^{\alpha}(z)$ is not suited for  constructive description of the $\alpha$-H{\"o}lder classes by means of complex polynomials.  This circumstance stimulated the introduction of  a modified scale $\rho_{1/n}^{\star\alpha}(z)$ \cite{Shirokov01,Andrievskii85,Andrievskii94}, which was used for constructive description of $\alpha$-H{\"o}lder classes in Jordan domains with non-empty interior.

In the case where the interior of $G$ is empty, i.e.,  if $ G=\Gamma$, the problem of a constructive description of H{\"o}lder (or H{\"o}lder-like) classes of functions defined on $\Gamma$ by means of their approximation by polynomials turned out to be more intricate. For example, if  $G=\Gamma_{\beta}\stackrel{def}{=} [-1,\, 0]\cup [0,\, e^{\beta}]$, $0 <\beta <\pi$, then a simple combination of $\rho_{1/n}(z)$ and $\rho_{1/n}^{\star}(z)$ cannot provide a constructive description of the $\alpha$-H{\"o}lder class \cite{Andrievskii94}. Even in the case of $\Gamma_{\beta}$, the answer is obtained with the help of a Cantor-like construction of a scale using both scales $\rho_{1/n}^{\alpha}(z)$ and $\rho_{1/n}^{\star\alpha}(z)$ \cite{Shirokov77}.

V.~V.~Andrievskii \cite{Andrievskii94} found an alternative approach to the problem of a constructive approximation of functional classes on   Jordan arcs. He used a uniform approximation of a function $f$ defined on a Jordan arc $L$ by polynomials $P_n$ along with uniform estimates of $P'_n(z)$ in a neighbourhood of $L$. We notice that  harmonic polynomials can also be used for a constructive description of H{\"o}lder-like classes of functions on continuums in $\mathbb{C}$ (V.~V.~Andrievskii, \cite{Andrievskii86,Andrievskii88}).

 We emphasize that all above-mentioned constructions of the scales $\rho_{1/n}^{\alpha}(z)$ and $\rho_{1/n}^{\star\alpha}(z)$ and constructive descriptions of H{\"o}lder classes on curves are applicable only for plane curve since each of these constructions uses  a conformal mapping of the complement $\mathbb{C}\setminus \overline{G}$ onto the exterior of the unit disc $\mathbb{D}$. However, the same problems can be considered for H{\"o}lder spaces on curves lying in arbitrary $\mathbb R^n$ or $\mathbb C^n$. 
 
 In the present paper, we obtain a constructive description of $H^{\alpha} (L)$ for   chord-arc curves $L$ lying in $\mathbb{R}^3$. 
As approximating functions, we use harmonic functions with certain estimates of their gradients in neighborhoods of a curve. The neighborhoods are connected with the rate of approximation -- they shrink when the approximation is getting better.  

The paper is organized as follows. In Section~2 we introduce notation and state our main results. Subsection~\ref{subsec:Theorem3}\, contains the proof of Theorem~3. Subsection~\ref{subsec:Theorem4}\, contains the proof of Theorem~4. Subsection~\ref{subsec:Theorem1}\, is concerned with the proof of the main result of the paper -- Theorem~1. Subsection~\ref{subsec:properties}\, is devoted to properties of a function $\upsilon_{2^{-n}}(M_0)$. Subsection~\ref{subsec:Theorem2}\, is devoted to the proof of Theorem~2.

{}
\section{Main results}

We say that a non-closed Jordan curve $L\subset \mathbb{R}^3$ has a chord-arc property (or is a chord-arc curve) if there exists a constant $C=C(L)$ such that the length of the subarc $L$ between $M_1$ and $M_2$ does not exceed $C\cdot\|M_1M_2\|$ for all points $M_1$, $M_2\in L$, $\|M_1M_2\|$ means the distance between $M_1$ and $M_2$ in $\mathbb{R}^3$.
We denote by $B_{r}(M)$ an open ball in $\mathbb{R}^3$ with center $M$ and   radius $r$ and put $\Omega_{\delta}(L)=\bigcup_{M\in L}  B_{\delta}(M)$. Let $H^{\omega}(L)$ be the space of all complex-valued functions $f$ that are defined on $L$ and satisfy  the condition $|f(M_2)-f(M_1)|\leq C_{f}\omega \big ( \|M_1M_2\|\big )$, where $\omega$ is a modulus of continuity with the property


\begin{equation}
\int_{0}^{x}\frac{\omega (t)}{t} dt \le C' \omega (x),\; x\int_{x}^{\infty}\frac{\omega (t)}{t^2} dt\le C'' \omega (x)
\label{f1}
\end{equation}
\noindent
(here and below we denote by $C$, $C'$, $C_1$, $\ldots$ various constants). One of our two  main results in the present paper is the following
theorem.

\begin{thm}
Assume that $L$ is a bounded non-closed chord-arc curve and $f\in H^{\omega}\big ( L \big )$. Then there exist constants $C_1=C_1 (f, L)$ and $C_2=C_2 (f, L)$ such that for every $\delta >0$ there exists a function $\upsilon_\delta$ harmonic in $\Omega_{\delta} (L)$ such that

\begin{equation}
\ |\upsilon_\delta(M)-f(M)|\leq C_1\omega(\delta), \; M\in L
\label{f2}
\end{equation}
\begin{equation}
\ |\nabla \upsilon_\delta(M)|\leq C_2 \frac{\omega(\delta)}{\delta}, \; M\in \Omega_{\delta}(L)\backslash\Omega_{\frac{\delta}{2}}(L)
\label{f3}
\end{equation}
\end{thm}
Theorem~1 may be called "a direct theorem" of approximation like many similar statements concerning approximation by polynomials, rational functions, etc. The "converse theorem" to   Theorem~1 is also valid: if we take a  unit vector $\vec{\ell}$, then (3) implies that

$$
|\upsilon_{\delta\ell}' (M)|\leq C_2\frac{\omega(\delta)}{\delta}, \, M\in \Omega_{\delta}(L)\backslash\Omega_{\frac{\delta}{2}}(L),\eqno(3')
$$

\noindent
and the maximum principle for a function $\upsilon_{\delta\ell}'$ harmonic in $\Omega_{\delta}(L)$ guarantees that   estimate (${3}^{\prime}$) is valid for $M\in\Omega_{\delta}$; this gives the estimate

$$
|\nabla \upsilon_{\delta}(M)|\leq C'_2 \frac{\omega(\delta)}{\delta}, \; M\in \Omega_{\delta}(L). \eqno(3'')
$$
Further, if $M_1, \, M_2\in L$ and $\|M_1M_2\| \leq \frac{\delta}{2}$, then the segment $\overrightarrow{M_{1}M_{2}}$ lies in $\Omega_{\delta}(L)$. Putting $\overline{\nu}=\frac{1}{\|M_1M_2\|}\cdot \; \overrightarrow{M_{1}M_{2}}$, we get

$$
f(M_2)-f(M_1)=\left( f(M_2)-\upsilon_{\delta}(M_2)\right)- \left(f(M_1)-\upsilon_{\delta}(M_1)\right)-
$$
$$
-\int\limits_0^1 \upsilon_{\delta\overline{\nu}}' (M_1 + \overline{\nu} \|M_1M_2\|t) dt. \eqno(3^{\circ})
$$

\noindent
So, if we suppose that a function $f$ can be approximated by functions $\upsilon_{\delta}$ as in $\eqref{f2}$ and $\eqref{f3}$, then ($3^{\circ}$) and ($3'$) imply that $f\in \mathbf{H}^{\omega}(L)$.
The constructive description of the space $\mathbf{H}^{\omega}(L)$ in terms of functions $\upsilon_{\delta}$ harmonic in  $\Omega_{\delta}(L)$ is in a sense strict. This is the assertion of the second main result.

\begin{thm}
Let $1>\delta_{k}>0$, $\delta_{k}\to 0$, $\delta_{k}$ be monotonically decreasing, $\ell_{k}\to +\infty$, and let the modulus of continuity $\omega (t)$ satisfy conditions $\eqref{f1}$. Then there exists a function $f_0 \in \mathbf{H}^{\omega}\left([A_0,\, B_0]\right)$, where $A_0=(-1,\, 0,\, 0)$ and $B_0=(1,\, 0,\, 0)$, that cannot be approximated  by functions $V_{k}$ harmonic in the domains $\Omega_{\ell_{k}\delta_{k}}\left([A_0,\, B_0]\right)$ in the following way:

\begin{equation}
\left|V_{k}(M)-f_0(M)\right|\leq C'_{1}\omega (\delta_{k}), \, M\in [A_0,\, B_0]
\end{equation}
\noindent
if the functions $V_{k}$ satisfy the condition

\begin{equation}
\left|\nabla V_{k}(M)\right|\leq C'_{2}\frac{\omega (\delta_{k})}{\delta_{k}}, \, M\in\Omega_{\ell_{k}\delta_{k}}\left([A_0,\, B_0]\right)
\end{equation}
\end{thm}

The proof of   Theorem~1 depends on a special type of an extension of a function $f$ from the curve $L$ to the entire space $\mathbb{R}^{3}$; we call this extension \emph{pseudoharmonic} by analogy with the widely-used pseudoanalytic extension due to E.~M.~Dyn'kin \cite{Dynkin77,Dynkin93}.

\begin{thm}
{\it Let  $f\in H^{\omega}\big ( L \big )$, where $\omega$ is a modulus of continuity satis\-fying assumption (1). Let  $O$ be the origin of $\mathbb{R}^{3}$. Then there is a function $f_0 \in C\left(\mathbb{R}^{3}\right)$ such that $f_0 |_{L}=f$, $f_0\in C^{2}\left(\mathbb{R}^{3}\setminus L\right)$, and}
\begin{equation}
\left|\nabla f_{0}(M)\right|=o(dist^{-1}(M,L)), \; o\; \text{is uniform on}\; \mathbb{R}^{3},
\end{equation}
\begin{equation}
f_{0}(M)\equiv 0, \;\text{for}\; \|\overrightarrow{OM}\| \geq R_{0}, \; \text{and}\; L\subset B_{R_{0}}(O)
\end{equation}
\begin{equation}
\left|\Delta f_{0}(M)\right|\leq C_0\frac{\omega\left( dist(M,L)\right)}{dist^2\left(M,L\right)}
\end{equation}
\end{thm}
In what follows, we call an extension $f_0$ of a function $f$ a pseudoharmonic extension of $f$.

\begin{thm}
{\it Assume that  a function $f\in C(L)$ has a pseudoharmonic extension satisfying conditions (6), (7), and (8). Then $f\in H^{\omega}\left( L \right)$}.
\end{thm}

Theorems~3 and 4 are exactly analogous to the theorems of E.~M.~Dyn'kin concerning pseudoanalytic extensions of   functions defined on domains in $\mathbb{C}$ \cite{Dynkin77,Dynkin93}.

\subsection{Proof of Theorem~3}
\label{subsec:Theorem3}

We begin with the proof of Theorem~3. Let $A$ be one of endpoints of the curve $L$ and let $B$ be the another one. In the sequel, we denote by $\ell (M_1,\, M_2)$   the length of the arc of $L$ with the endpoints $M_1$ and $M_2$. Let $\ell (A,\, B)=\Lambda$. We subdivide $L$ into $2^{n}$ arcs of equal length by the points $M_{kn}$, $0\leq k\leq 2^{n}$, $M_{0n}=A$, $M_{2^n,n}=B$, where the index  $k$ increases as the points $M_{kn}$ move in the direction from $A$ to $B$. The chord-arc property of $L$ implies the inequality

$$
\|\overrightarrow{M_{kn},\,M_{k+1,n}}\| \geq \frac{1}{C_0} \ell\left(M_{kn},\, M_{k+1,n}\right)=\frac{1}{C_0}\cdot 2^{-n} \Lambda \stackrel{\rm def}=\frac{1}{C_0} \Lambda_{n}.
$$

\noindent
We put

\begin{equation}
\Omega_{n}^\star\stackrel{\rm def}=\bigcup_{k=0}^{2^n} \overline{B}_{2\Lambda_{n}}(M_{kn}),
\end{equation}

\begin{equation}
\Omega_{n}\stackrel{\rm def}=\Omega_{n}^\star\setminus \overline{\Omega}_{n+1}^\star.
\end{equation}

\noindent
For $M\in \Omega_{n}$ we have the estimates

\begin{equation}
\frac{1}{2}\Lambda_{n}\leq dist\left(M,\,L\right)\leq 2\Lambda_{n}.
\end{equation}

\noindent
Let

\begin{equation}
\omega_{0n}=B_{2\Lambda_{n}}\left(M_{0n}\right)\cap{\Omega_{n}},
\end{equation}

\begin{equation}
\omega_{kn}=\left(B_{2\Lambda_{n}}\left(M_{kn}\right)\cap{\Omega_{n}}\right)\setminus \bigcup_{\nu=0}^{k-1} B_{2\Lambda_{n}}\left(M_{\nu n}\right),\; 1\leq k\leq 2^{n}
\end{equation}
($\omega_{kn}$ may be empty for some $k$ and $n$). We define the function $g$ as follows:

\begin{equation}
g(M)=\left\{
\begin{array}{rl}
f\left(M_{kn}\right),\, M\in\omega_{kn} \\
0, \, M\in \mathbb{R}^{3}\setminus\bigcup_{n=0}^{\infty} \Omega_{n}^{\star}
\end{array}\right.
\end{equation}

Let $d(M)=dist\left(M,\,L\right)$, $M\in \mathbb{R}^{3}\setminus L$ and $B_{\star}\left(M\right)=\overline{B}_{\frac{1}{8}d(M)}\left(M\right).$
\noindent
We need to control the distance $\|\overrightarrow{M_{kn},\, M_{k_1,\, n_1}}\|$ in the case where $M\in \omega_{kn}$, $M_1\in B_{\star}(M)\cap\omega_{k_1,\, n_1}$. We have $$2\Lambda_{n_1}\ge d(M_1)\ge d(M)-\|\overrightarrow{M M_1}\|\ge\frac{1}{2}\Lambda_{n}-\frac{1}{8} d(M)\ge \frac{1}{2}\Lambda_{n}-2\cdot\frac{1}{8}\Lambda_{n}=\frac{1}{4}\Lambda_{n},$$ from which we obtain $8\Lambda_{n_{1}}\ge\Lambda_{n}$, $-n_{1}+3\ge -n$, and $n_{1}\leq n+3$.

\noindent
Then we observe that
$$
\frac{1}{2}\Lambda_{n_1}\leq d(M_1)\leq d(M) +\|\overrightarrow{M M_1}\|\leq
$$
$$
\leq 2\Lambda_{n}+\frac{1}{8} d(M)\leq 2\Lambda_{n} + 2\cdot\frac{1}{8}\Lambda_{n} = 2\frac{1}{4}\Lambda_{n} <4\Lambda_{n} ,
$$

\noindent
hence $\Lambda_{n_1} < 8\Lambda_{n}$, $-n_1 < -n+3$, and $n_1\ge n-2$.

\noindent
Let $N, N_1\in L$  be such  that $\|\overrightarrow{MN}\| = d(M)$, $\|\overrightarrow{M_1 N_1}\| = d(M_1)$. Since $\|\overrightarrow{N M_{kn}}\| \leq 4\Lambda_{n}$, $\|\overrightarrow{N_1 M_{k_1 n_1}}\| \leq 4\Lambda_{n_1}$,
and
$$
\|N N_1\|\leq \|NM\|+\|M M_1\| + \|M_1 N_1\|\leq 2\Lambda_{n}+\frac{1}{8}\Lambda_{n} +2\Lambda_{n_1} \leq
$$
$$
\leq 2\Lambda_{n}+\frac{1}{8}\Lambda_{n} +2\cdot 8\Lambda_{n}=18\frac{1}{8}\Lambda_{n} < 19\Lambda_{n} ,
$$

\noindent
we have the estimates
\begin{equation}
\begin{split}
&\|M_{kn} M_{k_1 n_1}\|\leq \|M_{kn} N\| +\|N N_1\| +\|N_1 M_{k_1 n_1}\|\leq  \\
& \leq 4\Lambda_{n} + 19 \Lambda_{n} +4\Lambda_{n_1}\leq \left(4 + 19 +4\cdot 8\right) \Lambda_{n} < 55\Lambda_{n}
\end{split}
\end{equation}

\noindent
Inequality (15) and assumption ($8'$) imply the inequalities

\begin{equation}
\left| f\left(M_{kn}\right) - f\left(M_{k_1 n_1}\right)\right|\leq \omega\left(55\Lambda_{n}\right)\leq C\omega\left(\Lambda_{n}\right) .
\end{equation}
As a consequence of (16) and (14) we get the inequality
\begin{equation}
\left| g\left(M_1\right) - g\left(M\right)\right|\leq C\omega\left(d(M)\right).
\end{equation}
valid for all $M_{1}\in B_{\star}(M)$.
We define

\begin{equation}
g_{1}(M)=\frac{1}{\left| B_{\star}(M)\right|}\int\limits_{B_{\star}(M)} g\left( M_1\right) d{m_3}(M_1) ,
\end{equation}
where $\left| B_{\star}(M)\right|$  is the volume of the ball $B_{\star}(M)$ and $m_3$ is the 3-dimensional Lebesgue measure. Due to (18) and (17) we see that $g_1\in C\left(\mathbb{R}^{3}\setminus L\right)$ and

\begin{eqnarray}
&\left| g_1\left(M\right) - g\left(M\right)\right| =\left| g_1\left(M\right) - f\left(M_{kn}\right)\right|= \nonumber \\
&=\left|\frac{1}{\left| B_{\star}(M)\right|}\int\limits_{ B_{\star}(M)} g(M_1)\, d{m_3}\left( M_1\right) - \frac{1}{\left| B_{\star}(M)\right|}\int\limits_{ B_{\star}(M)} g(M)\, d{m_3}\left( M_1\right)\right|\leq \\
&\leq C\omega\left(d(M)\right) \nonumber .
\end{eqnarray}
The definition (14) and  estimate (19) imply that $g_{1}(M)\rightarrow f(M_{\star})$ as $M\rightarrow M_{\star}$, $M_{\star}\in L$. Hence  the function $g_1$ is continuous on $\mathbb{R}^{3}$ and vanishes outside a certain ball.

Now we construct   a characteristic $d_{0}(M)\approx d(M)$, but $d_{0}(M)$ is $C^2\left(\mathbb{R}^{3}\setminus L\right)$-smooth in contrast to $d(M)$, which is usually only Lip1 on $\mathbb{R}^{3}\setminus L$. Let $\sum_{n}^{} = \{ M\in \mathbb{R}^{3}\setminus L :\, 2^{n-1} < d(M)\leq 2^{n}\}$, $n\in \mathbb{Z}$. Since
$$
\left| d(M_2) - d(M_1)\right|\leq \|M_1 M_2\|,\; M_1, M_2 \in \mathbb{R}^{3}\setminus L ,
$$
\noindent
the balls $B_{r_1}\left( M_1\right)$ and $B_{r_2}\left( M_2\right)$ are disjoint if $r_1 < \frac{1}{4} d(M_1)$, $r_2 < \frac{1}{4} d(M_1)$, and $d(M_2)\ge 2d(M_1)$. Due to this observation,  the following functions are well defined:

\begin{equation}
d_{1}(M)=2^{n-1}, \; M\in\sum_{n}, \; n\in\mathbb{Z}
\end{equation}

\begin{equation}
d_{2}(M)=\frac{1}{\left| B_{\frac{1}{8}\cdot 2^{n-1}}(M)\right|}\int\limits_{ B_{\frac{1}{8}\cdot 2^{n-1}}(M)} d_{1}(\widetilde{M})\, dm_{3}(\widetilde{M}),
\end{equation}

\noindent
if $2^{n-1}\cdot\sqrt{2}<d(M)\leq \sqrt{2}\cdot 2^{n}=\frac{1}{\sqrt{2}}\cdot 2^{n+1}$.
We observe that definitions (20) and (21) imply the estimate $\|\operatorname{grad}d_{2}(M)\|\leq C$. Finally, we put

\begin{equation}
d_{0}(M)=\frac{1}{\left| B_{\frac{1}{8}\cdot 2^{n-1}}(M)\right|}\int\limits_{ B_{\frac{1}{8}\cdot 2^{n-1}}(M)} d_{2}(K)\, dm_{3}(K),
\end{equation}

\noindent
if $2^{n-1}\cdot\sqrt{2}<d(M)\leq \sqrt{2}\cdot 2^{n}$.

Equation (22) gives the required function $d_0$. We have the following estimates:

\begin{equation}
d_{0}(M)\asymp d(M), \; \|\operatorname{grad}d_{0}(M)\|\leq C
\end{equation}

\noindent
and

\begin{equation}
\|\operatorname{grad}^{2}d_{0}(M)\|\leq \frac{C}{d(M)},
\end{equation}

\noindent
which follow from  (22). Indeed, if $\bar{\lambda}$, $\bar{\mu}$ are arbitrary unit vectors, then (22) implies

$$
{d'_{0}}_{\bar{\lambda}}(M)=\frac{1}{\left| B_{\frac{1}{8}\cdot 2^{n-1}}(M)\right|}\int\limits_{B_{\frac{1}{8}\cdot 2^{n-1}}(M)} {d'_{2}}_{\bar{\lambda}}(K)\, dm_{3}(K),
$$

\noindent
which gives (23), and if $\bar{\nu}(K)$ is the outer unit normal to the sphere $S_{\frac{1}{8}\cdot 2^{n-1}}(M)$ at the point $K$,
then
\begin{equation}
{d''_{0}}_{\bar{\lambda}\bar{\mu}}(M)=\frac{1}{\left| B_{\frac{1}{8}\cdot 2^{n-1}}(M)\right|}\int\limits_{B_{\frac{1}{8}\cdot 2^{n-1}}(M)}\left(\bar{\mu}, \,\bar{\nu}(M) \right) {d'_{2}}_{\bar{\lambda}}(K)\, dS(K),
\end{equation}
where
$dS(K)$ denotes the Lebesgue measure on $S_{\frac{1}{8}\cdot 2^{n-1}}(K)$;   estimate (24) follows from (25). Let us notice that $d_{1}(M)\leq d(M)$, and, for $K\in B_{\frac{1}{8}\cdot 2^{n-1}}(M)$, we also have $d_{1}(K)\preccurlyeq d(M)$, hence $d_{2}(M)\preccurlyeq d(M)$. Moreover, (22) implies that $d_{0}(M)\preccurlyeq d(M)$. Finally, we define

\begin{equation}
g_{2}(M)=\frac{1}{\left| B_{\frac{1}{8} d_{0}(M)}(M)\right|}\int\limits_{B_{\frac{1}{8} d_{0}(M)}(M)} g_{1}(K)\, dm_{3}(K),
\end{equation}

\begin{equation}
g_{0}(M)=\frac{1}{\left| B_{\frac{1}{8} d_{0}(M)}(M)\right|}\int\limits_{B_{\frac{1}{8} d_{0}(M)}(M)} g_{2}(K)\, dm_{3}(K),
\end{equation}

\noindent
We notice  that definitions (20)--(22) imply the inequalities $$d_{1}(M)\ge\frac{1}{2} d(M),  \ d_{2}(M)\ge\frac{1}{2} d(M), \ d_{0}(M)\ge\frac{1}{2} d(M).$$
Let $B^{\star}(M)=B_{\frac{1}{8} d_{0}(M)}(M)$ and $r^{\star}(M)=\frac{1}{8} d_{0}(M)$. Using these estimates in the same  way as in (19), we get the estimates

\begin{equation}
\left| g_2(M)-g(M)\right|\leq C\omega (d(M))
\end{equation}

\noindent
and

\begin{equation}
\left| g_0(M)-g(M)\right|\leq C\omega (d(M)).
\end{equation}

\noindent
Let $\bar{\lambda}$ be a unit vector. We have
\begin{equation}
\begin{split}
g'_{2\bar{\lambda}}(M) &={{\left( g_2(N)-g(N) \right)}^{\prime}_{\bar{\lambda}\mid}}_{N=M}=  \\
&={{\left(\frac{1}{\left| B^{\star}(N)\right|}\int\limits_{B^{\star}(N)}\left( g_{1}(K)-g(M)\right)\, dm_{3}(K)\right)}'_{\bar{\lambda}\mid}}_{N=M}=  \\
&={{\left(\frac{1}{\left| B^{\star}(N)\right|}\right)}'_{\bar{\lambda}\mid}}_{N=M}\int\limits_{B^{\star}(M)} \left( g_1(K)-g(M)\right)\, dm_{3}(K) +  \\
&+\frac{1}{\left| B^{\star}(M)\right|} {{\left(\int\limits_{B^{\star}(N)}\left( g_{1}(K)-g(M)\right)\, dm_{3}(K)\right)}'_{\bar{\lambda}\mid}}_{N=M}=  \\
& =-\frac{{\left| B^{\star}(M)\right|}'_{\bar{\lambda}}}{{\left| B^{\star}(M)\right|}^{2}}\int\limits_{B^{\star}(M)}\left( g_{1}(K)-g(M)\right)\; dm_{3}(K) +  \\
& + \frac{1}{\left| B^{\star}(M)\right|}\int\limits_{\partial B^{\star}(M)} \left( \bar{n}(K),\bar{\lambda}\right)\left( g_1(K)-g(M)\right)\, dm_{2}(K),
\end{split}
\end{equation}

\noindent
where $\bar{n}(K)$ in the last integral is the unit vector of the outer normal to the sphere $\partial B^{\star}(M)$  and $dm_{2}(K)$ denotes the two-dimensional surface measure on the sphere $\partial B^{\star}(M)$.

Applying estimates (23) and (19) to   formula (30), we find that

\begin{equation}
\left| g'_{2\bar{\lambda}}(M)\right|\leq C\frac{\omega (d(M))}{d(M)},
\end{equation}

\noindent
hence

\begin{equation}
\left|\nabla g_{2}(M)\right|\leq C\frac{\omega (d(M))}{d(M)}.
\end{equation}
Repeating the same reasoning as in (30), we obtain by (28), (31), and (32) the following estimate  for $g_0$:
\begin{equation}
\left| g'_{0\bar{\lambda}}(M)\right|\leq C\frac{\omega (d(M))}{d(M)}.
\end{equation}

Let  $\bar{\lambda}$ and $\bar{\mu}$ be two arbitrary unit vectors. Then
\begin{eqnarray}
&{g''_{0}}_{\bar{\lambda}\bar{\mu}}(M)={{\left( g_0(N)-g(M) \right)}''_{\bar{\lambda}\bar{\mu}\mid}}_{N=M}= \nonumber \\
&={{\left(\frac{1}{\left| B^{\star}(N)\right|}\int\limits_{B^{\star}(N)}\left( g_{2}(K)-g(M)\right)\, dm_{3}(K)\right)}''_{\bar{\lambda}\bar{\mu}\mid}}_{N=M}  \nonumber \\
&={{\left(\frac{1}{\left| B^{\star}(N)\right|}\right)}''_{\bar{\lambda}\bar{\mu}\mid}}_{N=M}\int\limits_{B^{\star}(M)} \left( g_2(K)-g(M)\right)\, dm_{3}(K)  \nonumber \\
&+{{\left(\frac{1}{\left| B^{\star}(N)\right|}\right)}'_{\bar{\lambda}\mid}}_{N=M}
{{\left(\int\limits_{B^{\star}(N)} \left( g_2(K)-g(M)\right)\, dm_{3}(K)\right)}'_{\bar{\mu}\mid}}_{N=M} \nonumber \\
&+{{\left(\frac{1}{\left| B^{\star}(N)\right|}\right)}'_{\bar{\mu}\mid}}_{N=M}
{{\left(\int\limits_{B^{\star}(N)} \left( g_2(K)-g(M)\right)\, dm_{3}(K)\right)}'_{\bar{\lambda}\mid}}_{N=M} \nonumber \\
&+\frac{1}{\left| B^{\star}(M)\right|} {{\left(\int\limits_{B^{\star}(N)}\left( g_{2}(K)-g(M)\right)\, dm_{3}(K)\right)}''_{\bar{\lambda}\bar{\mu}\mid}}_{N=M}  \\
&=-{{\left(\frac{{\left| B^{\star}(N)\right|}'_{\bar{\lambda}}}{{\left| B^{\star}(N)\right|}^{2}}\right)}'_{\bar{\mu}\mid}}_{N=M}\int\limits_{B^{\star}(M)}\left( g_{2}(K)-g(M)\right)\; dm_{3}(K) \nonumber \\
&-\frac{{{\left| B^{\star}(N)\right|}'_{\bar{\lambda}\mid}}_{N=M}}{{\left| B^{\star}(M)\right|}^{2}} \nonumber
     \int\limits_{\partial B^{\star}(M)}\left(\left( \bar{n}(K),\bar{\mu}\right)+\left( r^{\star}(M)\right)'_{\bar{\mu}}\right)\cdot\left( g_2(K)-g(M)\right)\, dm_{2}(K) \nonumber  \\
&-\frac{{{\left| B^{\star}(N)\right|}'_{\bar{\mu}\mid}}_{N=M}}{{\left| B^{\star}(M)\right|}^{2}}
     \int\limits_{\partial B^{\star}(M)}\left(\left( \bar{n}(K),\bar{\lambda}\right)+\left( r^{\star}(M)\right)'_{\bar{\lambda}}\right)\cdot\left( g_2(K)-g(M)\right)\, dm_{2}(K) \nonumber \\
&+ \frac{1}{\left| B^{\star}(M)\right|}
     {\left( \int\limits_{\partial B^{\star}(N)}\left(\left( \bar{n}(K),\bar{\lambda}\right)+\left( r^{\star}(N)\right)'_{\bar{\lambda}}\right)
		\left( g_2(K)-g(M)\right)\, dm_{2}(K)\right)}'_{{\left.\bar{\mu}\right|}_{N=M}} \nonumber
\end{eqnarray}

Now we take into account that

\begin{equation}
\begin{split}
&{{\left( \int\limits_{\partial B^{\star}(N)}\left(\left( \bar{n}(K),\bar{\lambda}\right)+\left( r^{\star}(N)\right)'_{\bar{\lambda}}\right)\left( g_2(K)-g(M)\right)\, dm_{2}(K)\right)}'_{\bar{\mu}\mid}}_{N=M}= \\
&=\int\limits_{\partial B^{\star}(M)}{{\left( {r^{\star}(N)}'_{\bar{\lambda}}\right)}'_{\bar{\mu}\mid}}_{N=M}\cdot\left( g_2(K)-g(M)\right)\, dm_{2}(K)+ \\
& + 2\int\limits_{\partial B^{\star}(M)}\frac{{( r^{\star}(M) )}'_{\bar{\lambda}}{( r^{\star}(M) )}'_{\bar{\mu}}}{r^{\star}(M)}\left( g_2(K)-g(M)\right)\, dm_{2}(K)+  \\
&+\int\limits_{\partial B^{\star}(M)}{(r^{\star}(N))}'_{\bar{\lambda}}\cdot {(g_2(K))}'_{\bar{\mu}}\, dm_{2}(K)+ \\
&+2\int\limits_{\partial B^{\star}(M)} \frac{\left(\bar{n} (K),\bar{\lambda}\right){(r^{\star}(M))}'_{\bar{\mu}}}{r^{\star}(M)}\left( g_2(K)-g(M)\right)\, dm_{2}(K)+  \\
&+\int\limits_{\partial B^{\star}(M)}\left(\left( \bar{n}(K),\bar{\lambda}\right)+{\left( r^{\star}(M)\right)}'_{\bar{\lambda}}\right) {(r^{\star}(M))}'_{\bar{\mu}} {(g_2(K))}'_{\bar{n}(K)}\, dm_{2}(K).
\end{split}
\end{equation}
Combining   estimates (23), (24), (28), and (33) and equalities (34) and (35), we find that

\begin{equation*}
\left| {g_0}_{\bar \lambda\bar \mu }(M)\right|\leq C\frac{\omega (d(M))}{d^2(M)},
\end{equation*}

\noindent
which implies

\begin{equation*}
\left|\nabla^2 g_{0}(M)\right|\leq C\frac{\omega (d(M))}{d^2(M)},
\end{equation*}

\noindent
and finally,

\begin{equation}
\left|\Delta g_{0}(M)\right|\leq C\frac{\omega (d(M))}{d^2(M)}.
\end{equation}

Inequalities (29), (33), and (36) conclude the proof of   Theorem~3 with a slight change in notation: we have produced a required function $g_0$.

\subsection{Proof of Theorem~4}
\label{subsec:Theorem4}


Now we  proceed to the proof of Theorem~4. Consider the sets $\Omega_{\!n}^\star$ and $\Omega_{\!n}$ defined in (9) and (10). The boundaries of $\Omega_n$ and $\Omega_{\!n}^\star$ consist of a finite number of subsets of spheres of radii $2\Lambda_{\!n}$ and $\Lambda_{\!n}$; the total area of these spheres is

\begin{equation}
4\pi \left( (2^{n+1}+1)\cdot\Lambda_{\!n}^{2} + (2^n+1)\cdot 4\Lambda_{\!n}^{2}\right)\leq C\cdot 2^n\cdot ({2^{-n}})^2 = C\cdot 2^{-n}
\end{equation}
We fix a point $M_0\in\mathbb{R}^3\setminus L$ and choose $n$ such that $M_0\notin\Omega_{\!n}^\star$. Assume that $f_0$ is a pseudoharmonic extension of $f$ and that  $R_0$ is chosen so large that $f_0(M)\equiv 0$ outside  the ball $B_{R_0}({O})$ and $M_0\in B_{R_0}({O})$. We denote by $\Sigma_{\!n}$ the connected component of the set $B_{R_0}({O})\setminus \Omega_{\!n}$ containing the point $M_0$. Now we use   the classical formula

\begin{equation}
\begin{split}
&f_0(M_0)=\frac{1}{4\pi} \int\limits_{\partial\sum_{n}} \left(f_0(M)\right)'_{\bar{n} (M)} \, \frac{1}{\rho_{M_0}(M)} \, dS(M)- \\
&-\frac{1}{4\pi}\int\limits_{\partial\sum_{n}} \, f_0(M)\left(\frac{1}{\rho_{M_0}(M)}\right)'_{\bar{n} (M)} \, dS(M)- \\
&-\frac{1}{4\pi}\int\limits_{\partial\sum_{n}} \,\frac{\Delta f_0(M)}{\rho_{M_0}(M)} \, dm_{3}(M),
\end{split}
\end{equation}

\noindent
where $\rho_{\!M_0}(M)\stackrel{\rm def}= \|M_0 M\|$ $\vec n(M)$ is the outer unit normal at $M\in\partial \Sigma_{\!n}$  to the domain $\Sigma_{\!n}$, $dS(M)$ is the two-dimensional measure on $\partial \Sigma_{\!n}$,  and $m_3$ is the three-dimensional Lebesgue measure in $\mathbb{R}^{\!3}$.

We take into account that $f_0(M)\equiv 0$ and ${\left( f_0(M)\right)}'_{\pi (M)}\equiv 0$ for $M\in \partial B_{R_0}({O})$. This implies that the integrals in (38) are calculated over the domain $\partial \Sigma_{\!n}\bigcap\partial \Omega_{\!n}^\star$ whose two-dimensional measure does not exceed $c\cdot 2^{-n}$. The    construction of $\Omega_{\!n}^\star$ gives the estimates $c'\cdot 2^{-n}\leq d(M)\leq c''\cdot 2^{-n}$, $M\in \partial \Omega_{\!n}^\star$, with some constants $c',\; c'' > 0$, and   condition (6) yields a sequence $\{\alpha_n\}_{n=1}^{\infty}$, $\alpha_n\to 0$, such that


\begin{equation}
\left|\ \left( f_0(M)\right)'_{\bar{n} (M)}\right|\leq C \alpha_{n}\left( d(M)\right)^{-1}, \, M\in \Sigma_n.
\end{equation}
Using (39) and the above argument, we obtain


\begin{equation}
\left|\frac{1}{4\pi} \int\limits_{\partial\sum_{n}} \left(f_0(M)\right)'_{\bar{n} (M)} \, \frac{1}{\rho_{M_0}(M)} \, dS(M) \right| \leq C\alpha_{n} \cdot 2^{n}\cdot 2^{-n}=C\alpha_{n}\\
\end{equation}

\noindent
and

\begin{equation}
\left|-\frac{1}{4\pi} \int\limits_{\partial\sum_{n}} f_0(M) \cdot \, \left(\frac{1}{\rho_{M_0}(M)}\right)'_{\bar{n} (M)} \, dS(M) \right| \leq C\cdot 2^{-n}. \end{equation}
Formula (38) and estimates (40) and (41) imply the relation
\begin{equation}
\begin{split}
f_0(M_0)=-\frac{1}{4\pi} \int\limits_{\sum_{n}}  \frac{\Delta f_0(M)}{\rho_{M_0}(M)} \, dm_3(M) + O(\alpha_{n} + 2^{-n})
\end{split}
\end{equation}
Passing to the limit in (42), we get
\begin{equation}
\begin{split}
f_0(M_0)=-\frac{1}{4\pi} \int\limits_{B_{R_0}(\mathbb{O})}  \frac{\Delta f_0(M)}{\rho_{M_0}(M)} \, dm_3(M).
\end{split}
\end{equation}

We will check below that the integral in   (43) is continuous on $\mathbb{R}^{3}$. Equality (43) and the continuity of both sides of it on $\mathbb{R}^{3}$ allows us to take in (43) an arbitrary point $M_0$ of $\mathbb{R}^{3}$. In particular, we can take $M_0\in L$. Bearing this in mind, we take $M_1, M_2 \in L$, $M_1\neq M_2$ and obtain

\begin{equation}
\begin{split}
&f(M_2)-f(M_1)=\frac{1}{4\pi} \int\limits_{B_{R_0}(\mathbb{O})}  \frac{\Delta f_0(M)}{\rho_{M_1}(M)} \, dm_3(M)-\\
&- \frac{1}{4\pi} \int\limits_{B_{R_0}(\mathbb{O})}  \frac{\Delta f_0(M)}{\rho_{M_2}(M)} \, dm_3(M) =\\
&= \frac{1}{4\pi} \int\limits_{B_{2\|M_1M_2\|}(M_1)}  \frac{\Delta f_0(M)}{\rho_{M_1}(M)} \, dm_3(M) -\\
&-\frac{1}{4\pi} \int\limits_{B_{2\|M_1M_2\|}(M_1)}  \frac{\Delta f_0(M)}{\rho_{M_2}(M)} \, dm_3(M)+\\
&+\frac{1}{4\pi} \int\limits_{B_{R_0}(\mathbb{O})\setminus B_{2\|M_1M_2\|}(M_1)}  \left(\frac{1}{\rho_{M_1}(M)}-\frac{1}{\rho_{M_2}(M)}\right) \Delta f_0(M)\, dm_3(M) \\
& \stackrel{\rm def}= I_1-I_2+I_3.
\end{split}
\end{equation}

We remind that we assume relations (6), (7), and (8). Using them we get

\begin{equation}
\begin{split}
&|I_2|\leq \frac{1}{4\pi} \int\limits_{B_{3\|M_1M_2\|}(M_2)}  \frac{|\Delta f_0(M)|}{\rho_{M_2}(M)} \, dm_3(M)\leq \\
&C\int\limits_{B_{3\|M_1M_2\|}(M_2)}  \frac{\omega(d(M))}{d^2(M) \rho_{M_2}(M)} \, dm_3(M)=\\
&=C\sum_{n=0}^{\infty}\int\limits_{B_{3\cdot 2^{-n}\|M_1M_2\|}(M_2)\setminus B_{3\cdot 2^{-n-1}\|M_1M_2\|}(M_2)}  \frac{\omega(d(M))}{d^2(M) \rho_{M_2}(M)} \, dm_3(M)\leq \\
&C\sum_{n=0}^{\infty}\frac{2^n}{\|M_1M_2\|}\int\limits_{B_{3\cdot 2^{-n}\|M_1M_2\|}(M_2)\setminus B_{3\cdot 2^{-n-1}\|M_1M_2\|}(M_2)}  \frac{\omega(d(M))}{d^2(M)} \, dm_3(M)\leq \\
&C\sum_{n=0}^{\infty}\frac{2^n}{\|M_1M_2\|}\int\limits_{B_{3\cdot 2^{-n}\|M_1M_2\|}(M_2)}  \frac{\omega(d(M))}{d^2(M)} \, dm_3(M)
\end{split}
\end{equation}

Without loss of generality, we may assume that $\|M_1M_2\|\leq \widetilde{C}\|AB\|$ with a constant $\widetilde{C}$ such that   $B_{3\|M_1M_2\|}(M_2)\subset \Omega_{0}^\star$ for  $M_2\in L$, where $\Omega_{0}^\star$ is the set defined in (9).
Let $\sigma_{nk}=B_{3\cdot 2^{-n}\|M_1M_2\|} (M_2)\bigcap \Omega_{k}$, where the sets $\Omega_{k}$ are defined in (10). Then we can rewrite  a summand in (45) in the following way:

\begin{equation}
\begin{split}
&\int\limits_{B_{3\cdot 2^{-n}\|M_1M_2\|}(M_2)}  \frac{\omega(d(M))}{d^2(M)} \, dm_3(M)= \\
&=\sum_{k=0}^{\infty}\int\limits_{\sigma_{nk}} \frac{\omega(d(M))}{d^2(M)} \, dm_3(M) = \\
&=\sum_{k=k(n)}^{\infty}\int\limits_{\sigma_{kn}}  \frac{\omega(d(M))}{d^2(M)} \, dm_3(M).
\end{split}
\end{equation}

The index $k(n)$ in  (46) means the smallest $k$ such that $\Omega_k\cap B_{3\cdot 2^{-n}\|M_1M_2\|}(M_2)\neq \varnothing$. Inequalities (11) imply the following  important estimates:

\begin{equation}
2^{-k(n)}\asymp 2^{-n}\cdot \|M_1M_2\|,
\end{equation}

\begin{equation}
d(M)\asymp 2^{-k}, \, M\in \sigma_{kn}
\end{equation}


Let $\widetilde{\sigma}_{nk}=B_{3\cdot 2^{-n}\|M_1M_2\|} (M_2)\bigcap \Omega_{k}^\star$, then $\sigma_{nk}\subset \widetilde{\sigma}_{nk}$ and $m_3 \sigma_{nk}\leq m_3 \widetilde{\sigma}_{nk}$. Since $\Omega_{\nu}\bigcap B_{3\cdot 2^{-n}\|M_1M_2\|} (M_2)=\varnothing$ for $\nu < k(n)$, we see that, for $k\ge k(n)$, the center of each ball constituent  of $\widetilde{\sigma}_{nk}$ of radius $2^{-k}$  lies on a subarc of $L$ of  length $\leq C\cdot 2^{-n}\cdot \|M_1M_2\|$, which implies that the number $N_{n,k}$ of such balls does not exceed $C\cdot 2^{-n}\cdot \|M_1M_2\|\cdot 2^k$. Hence
\begin{equation}
m_3 \widetilde{\sigma}_{nk}\leq C N_{nk}\cdot 2^{-3k}\leq C\cdot 2^{-n}\cdot\|M_1M_2\|\cdot 2^k \cdot 2^{-3k}=C\cdot 2^{-n-2k}\|M_1M_2\|.
\end{equation}
Finally, combining   estimates (47), (48), and (49), we obtain

\begin{equation}
\begin{split}
&\sum_{k=k(n)}^{\infty}\int\limits_{\sigma_{kn}} \frac{\omega(d(M))}{d^2(M)} \, dm_3(M)\leq  \\
&C\sum_{k=k(n)}^{\infty} 2^{2k}\omega(2^{-k}) m_3(\widetilde{\sigma}_{nk})\leq \\
&C\sum_{k=k(n)}^{\infty} 2^{2k}\omega(2^{-k})\cdot 2^{-n-2k}\|M_2M_2\|=\\
&=C 2^{-n}\|M_1M_2\|\cdot\sum_{k=k(n)}^{\infty} \omega(2^{-k}).
\end{split}
\end{equation}

The first assumption in (1) concerning $\omega(t)$ gives the inequality

\begin{equation}
\sum_{k=k(n)}^{\infty} \omega(2^{-k})\leq C\omega(2^{-k(n)})\leq C\omega(2^{-n}\cdot \|M_1M_2\|).
\end{equation}

So, formulas (46), (50), and (51) imply the estimate
\begin{equation}
\int\limits_{B_{3\cdot 2^{-n}\|M_1M_2\|} (M_2)} \frac{\omega(d(M))}{d^2(M)} \, dm_3(M)\leq C 2^{-n}\cdot\|M_1M_2\|\cdot\omega(2^{-n}\cdot \|M_1M_2\|).
\end{equation}

Let us substitute (52) into (45). Using (1), we obtain

\begin{equation}
\begin{split}
&\sum_{n=0}^{\infty}\frac{2^n}{\|M_1M_2\|}\int\limits_{B_{3\cdot 2^{-n}\|M_1M_2\|} (M_2)} \frac{\omega(d(M))}{d^2(M)} \, dm_3(M)\leq \\ &C\sum_{n=0}^{\infty}\frac{2^n}{\|M_1M_2\|}\cdot 2^{-n}\cdot\|M_1M_2\|\omega(2^{-n}\cdot\|M_1M_2\|)= \\
&=C\sum_{n=0}^{\infty} \omega(2^{-n}\cdot\|M_1M_2\|)\leq C\omega(\|M_1M_2\|),\\
\end{split}
\end{equation}
which means that $|I_2|\leq C\omega(\|M_1M_2\|)$.

The same arguments show that $|I_1|\leq C\omega(\|M_1M_2\|)$. To estimate the term $I_3$, we  use the second part of   assumption (1) concerning the function $\omega(t)$. We notice that, for all $M\notin B_{2\|M_1M_2\|(M_1)}$,
we have the inequality
\begin{equation}
\left |\frac{1}{\rho_{M_1}(M)}-\frac{1}{\rho_{M_2}(M)}\right |\leq C\frac{\|M_1M_2\|}{\rho_{M_1}^2(M)}.
\end{equation}
Now, using (54) and (8), we obtain

\begin{equation}
\begin{split}
&|I_3|\leq C\int\limits_{B_{R_0}(\mathbb{O})\setminus B_{2\|M_1M_2\|}(M_1)}  \frac{\|M_1M_2\|}{\rho_{M_1}^2(M)} |\Delta f_0(M)|\, dm_3(M)\leq \\
&C\int\limits_{B_{R_0}(\mathbb{O})\setminus B_{2\|M_1M_2\|}(M_1)}  \frac{\|M_1M_2\|}{\rho_{M_1}^2(M)} \frac{\omega(d(M))}{d^2(M)}\, dm_3(M)\leq \\
&C\sum_{n=1}^{\infty}\int\limits_{(B_{2^{n+1}\|M_1M_2\|}(M_1)\setminus B_{2^{n}\|M_1M_2\|}(M_1))\bigcap B_{R_0}(\mathbb{O})}
\frac{\|M_1M_2\|}{\rho_{M_1}^2(M)} \frac{\omega(d(M))}{d^2(M)}\, dm_3(M)\leq \\
&C\sum_{n=1}^{\infty}\frac{\|M_1M_2\|}{2^{2n}\|M_1M_2\|^{2}}\int\limits_{(B_{2^{n+1}\|M_1M_2\|}(M_1)\setminus B_{2^{n}\|M_1M_2\|}(M_1))\bigcap B_{R_0}(\mathbb{O})}\frac{\omega(d(M))}{d^2(M)}\, dm_3(M)\stackrel{\rm def}=\\
&C\frac{1}{\|M_1M_2\|}\sum_{n=1}^{\infty}\frac{1}{2^{2n}} C_n.
\end{split}
\end{equation}

Now, repeating the same reasoning as we used to get (47)--(52), we obtain the estimate
\begin{equation}
C_n\leq C\cdot 2^{n}\cdot \|M_1M_2\|\cdot \omega(2^{n}\|M_1M_2\|).
\end{equation}

Combining (55) and (56), we see that

\begin{equation}
\begin{split}
&|I_3|\leq C\frac{1}{\|M_1M_2\|}\sum_{n=1}^{\infty}\frac{1}{2^{2n}}\cdot 2^{n}\cdot \|M_1M_2\|\cdot \omega(2^{n}\|M_1M_2\|)=\\
&=C\sum_{n=1}^{\infty}\frac{\omega(2^{n}\|M_1M_2\|)}{2^n}\leq C\omega(\|M_1M_2\|).
\end{split}
\end{equation}
We made use of the second part of   condition (1) in the last inequality in (57). So, we have proved that $|I_1|,\,|I_2|,\,|I_3|\leq C\omega(\|M_1M_2\|)$, which together with Proposition (44) finishes the proof of  Theorem~4.

\subsection{Proof of Theorem~1}
\label{subsec:Theorem1}


We start with some geometrical observations. We divide $L$ by the points $$A=M_{0n},\,M_{1n},\,\ldots M_{2^n n}=B$$ as we did in the definitions (9) and (10) of the domains $\Omega_{n}^{\star}$ and $\Omega_{n}$. Let $\Lambda_{n}=2^{-n}\cdot |\Lambda|$,   $C_1\ge 1,\, 0\leq k_0\leq 2^n$, and $B[C_1]=B_{C_1\cdot\Lambda_{n}}(M_{k_0 n})$. Let  $P_0,\, P_1\in \partial B[C_1]\bigcap L$ be such that the subarc $L(P_0,\, P_1)$ of $L$ with the endpoints $P_0$ and $P_1$ is the biggest one if $\partial B[C_1]\bigcap L$ contains more than two points. Then we have  $L(P_0,P_1)\leq C_0\cdot 2C_1\cdot \Lambda_{n}$, and there are at most $[2C_0C_1]+2\leq 2(C_0+1)C_1$ subarcs $L(M_{k,2^n},\, M_{k+1,2^n})$ intersecting $L(P_0,P_1)$. Then it is clear that $$m_3(B[C_1]\bigcap \Omega_{n}^{\star})\leq 2(C_0+1)C_1\cdot\frac{4}{3}\pi\cdot (2\Lambda_{n})^3=2(C_0+1)C_1\cdot \frac{32}{3}\pi\Lambda_{n}^{3}.$$

The volume of $B[C_1]$ is equal to $\frac{4}{3}\pi\cdot C_1^{3}\cdot\Lambda_{n}^3$. Therefore,  we can choose  $C_1$ such that $m_3(B[C_1])\ge 2m_3(B[C_1]\bigcap\Omega_{n}^{\star})$. We introduce the sets $\beta_{kn},\, 0\leq k\leq 2^n$, as follows: $\beta_{0n}=B_{2\Lambda_{n-2}(M_{0n})},\, \beta_{kn}=B_{2\Lambda_{n-2}}(M_{kn})\setminus\bigcup_{\nu=1}^{k-1}B_{2\Lambda_{n-2}}(M_{\nu n})$. We take a constant $C_1$ in such a way that $m_3(B_{C_1\Lambda_{n}}(M_{kn})\setminus\Omega_{n-2}^{\star})\ge \frac{1}{2}m_3(B_{C_1\Lambda_{n}}({O}))$.

The above arguments show that we can choose   $C_1$ depending only on $C_0$. Due to estimates (11) we obtain that the inequality $d(M)\ge 2^{-n+1}$ is  valid for all $M\in B_{C_1(\Lambda_{0})}(M_{kn})\setminus\Omega_{n-2}^{\star}$. On the other hand, $d(M)\leq C_1 2^{-n}|\Lambda|$.

Now we proceed   to the definition of   $\upsilon_{2^{-n}}(M)$. Using (46)--(52),
 we obtain

\begin{equation}
\int\limits_{B_{2\Lambda_{n-2}}(M_{kn})}\frac{\omega(d(M))}{d^2(M)}\, dm_3(M)\leq C\Lambda_{n-2}\omega(\Lambda_{n-2})
\end{equation}

Inequality (58) and the definition of the set $\beta_{kn}\subset B_{2\Lambda{n-2}}(M_{kn})$ imply

\begin{equation}
\int\limits_{\beta_{kn}}\frac{\omega(d(M))}{d^2(M)}\, dm_3(M)\leq C\Lambda_{n-2}\omega(\Lambda_{n-2}).
\end{equation}

Now we  apply Theorem~4 and construct a pseudoharmonic extension $f_0(M)$ of   $f$. Then (8) and (59) give the relation

\begin{equation}
\int\limits_{\beta_{kn}}\, \Delta f_0(M)\, dm_3 = C_{kn}\Lambda_{n-2}\omega(\Lambda_{n-2}),
\end{equation}

\noindent
where $|C_{kn}|\leq C$ for all $n$ and $k$, $0\leq k\leq 2^{n-2}$. We denote by $\chi_{kn}$ the characteristic function of the set $B_{C_1\Lambda_{n}}(M_{kn})\setminus \Omega_{n-2}^{\star}$ and put

\begin{equation}
\phi_{kn}(M)=\gamma_{kn}\, \Lambda^{-2}_{n}\chi_{kn}(M)\omega (\Lambda_{n}),
\end{equation}

\noindent
where $\gamma_{kn}$ satisfies the condition

\begin{equation}
\int\limits_{\beta_{kn}}\, \Delta f_0(M)\, dm_3 + \int\limits_{\mathbb{R}^3}\,\phi_{kn}(M)\, dm_3(M)=0.
\end{equation}
Taking into account (60) and (61) and the definition of the constant $C_1$, we obtain that $|\gamma_{kn}|\leq C$, where $C$ is independent of $k$ and $n$. Further, we define

\begin{equation}
\Phi_{n}=\sum_{k=0}^{2^{n-2}}\, \phi_{kn}(M).
\end{equation}
Preserving the notation $\rho_{M_0}(M)=\|M_0 M\|$, we define  the function $\upsilon_{2^{-n}}(M_0)$ as follows:
\begin{equation}
\begin{split}
&\upsilon_{2^{-n}}(M_0)=-\frac{1}{4\pi}\int\limits_{\mathbb{R}^3\setminus\Omega_{n-2}^{\star}}\,\frac{\Delta f_0(M)}{\rho_{M_0}(M)}\, dm_3(M)+ \\
&+\frac{1}{4\pi}\int\limits_{\mathbb{R}^3}\,\frac{\Phi_{n}(M)}{\rho_{M_0}(M)}\, dm_3(M)
\end{split}
\end{equation}

\subsection{Properties of a function $\upsilon_{2^{-n}}(M_0)$}
\label{subsec:properties}


Inequality (11) applied to the set $\Omega_{n-2}^{\star}$ shows that $(\SP \, \Delta f_0)\cap\Omega_{2^{-n + 1}}(L)=\varnothing$ and $(\SP \,  \Phi_n)\cap\Omega_{2^{-n+1}}(L)=\varnothing$. By (64) the function $\upsilon_{2^{-n}}$ is harmonic in $\Omega_{2^{-n+1}}(L)$. Assume that $M_0\in L$. Then, using (43) and (64), we get

\begin{equation}
\begin{split}
&\upsilon_{2^{-n}}(M_0)-f(M_0)=-\frac{1}{4\pi}\int\limits_{\mathbb{R}^3\setminus\Omega_{n-2}^{\star}}\,\frac{\Delta f_0(M)}{\rho_{M_0}(M)}\, dm_3(M)+ \\
&+\frac{1}{4\pi}\int\limits_{\mathbb{R}^3}\,\frac{\Phi_{n}(M)}{\rho_{M_0}(M)}\, dm_3(M)+\frac{1}{4\pi}\int\limits_{\mathbb{R}^3}\,\frac{\Delta f_0(M)}{\rho_{M_0}(M)}\, dm_3(M)=\\
&=\frac{1}{4\pi}\int\limits_{\Omega_{n-2}^{\star}}\,\frac{\Delta f_0(M)}{\rho_{M_0}(M)}\, dm_3(M)+\frac{1}{4\pi}\int\limits_{\mathbb{R}^3}\,\frac{\Phi_{n}(M)}{\rho_{M_0}(M)}\, dm_3(M)=\\
&=\sum_{k=0}^{2^{n-2}}\,\left(\frac{1}{4\pi}\int\limits_{\beta_{kn}}\,\frac{\Delta f_0(M)}{\rho_{M_0}(M)}\, dm_3(M)+\frac{1}{4\pi}\int\limits_{\mathbb{R}^3}\,\frac{\phi_{kn}(M)}{\rho_{M_0}(M)}\, dm_3(M)\right).
\end{split}
\end{equation}

Let $M_0$ belong  to the closed subarc $L\left (M_{k_0,n-2},\, M_{k_{0}+1,n-2} \right )$ of $L$ with the endpoints $M_{k_0,n-2}$ and $M_{k_{0}+1,n-2}$. By (62), we get
\begin{equation}
\sum_{k=0}^{2^{n-2}}=\sum_{k=0}^{k_{0}-2}+\sum_{k=k_{0}-1}^{k_{0}+2}+\sum_{k=k_{0}+3}^{2^{n-2}}\stackrel{\rm def}=\Sigma_1+\Sigma_2+\Sigma_3.
\end{equation}

Now in the same way as in (46)--(53), we get the estimates
\begin{equation}
\begin{split}
&\left|\frac{1}{4\pi}\int\limits_{\beta_{kn}}\,\frac{\Delta f_0(M)}{\rho_{M_0}(M)}\, dm_3(M)\right| \\
&\leq C \int\limits_{B_{2\Lambda_{n-2}(M_{k,n-2})}}\,\frac{\omega(d(M))}{d^{2}(M)\cdot 2^{-n+1}}\, dm_3(M) \\
&\leq C\omega(2^{-n+2})\leq C\omega(2^{-n})
\end{split}
\end{equation}

\noindent
for $k_0-1\leq k\leq k_0+2$ because $\rho_{M_0}(M)\ge 2^{-n+1}$ for $M_0\in L$ and $M\in\beta_{kn}$.
Moreover, for all $k$, $0\leq k\leq 2^{n-2}$, we have the inequalities
\begin{equation}
\begin{split}
&\left|\frac{1}{4\pi}\int\limits_{\mathbb{R}^3}\,\frac{\phi_{kn}(M)}{\rho_{M_0}(M)}\, dm_3(M)\right|=\frac{|\gamma_n|}{4\pi\, \Lambda^{2}_{n}}\omega(\Lambda_n)\int\limits_{B_{C_1\Lambda_n}(M_{kn})\setminus\widetilde{\Omega}_{n-2}}\, \frac{dm_3(M)}{\rho_{M_0}(M)}\leq \\
&\leq \frac{|\gamma_{kn}|}{4\pi\, \Lambda^{2}_{n}} \omega(\Lambda_n)\cdot \frac{1}{2^{-n+1}}\cdot m_3\left( B_{C_1\Lambda_n}(M_{kn})\setminus\widetilde{\Omega}_{n-2}\right)\leq\\
&\leq C\frac{\omega(\Lambda_n)}{\Lambda^{3}_{n}}\cdot \Lambda^{3}_{n}\leq C\omega(2^{-n}).
\end{split}
\end{equation}
\noindent
Relations  (67) and (68) imply that
\begin{equation}
|\Sigma_2|\leq C\omega(2^{-n}).
\end{equation}
Let us suppose now that $k\leq k_0-2$ or $k\ge k_0+3$. Then we transform the summands in $\Sigma_1$ or $\Sigma_3$ as follows:


\begin{equation}
\begin{split}
&\frac{1}{4\pi}\int\limits_{\beta_{kn}}\,\frac{\Delta f_0(M)}{\rho_{M_0}(M)}\, dm_3(M)+\frac{1}{4\pi}\int\limits_{\mathbb{R}^3}\frac{\phi_{kn}(M)}{\rho_{M_0}(M)}\, dm_3(M) = \\
&=\frac{1}{4\pi}\int\limits_{\beta_{kn}}\frac{\Delta f_0(M)}{\rho_{M_0}(M_{k,n-2})}\, dm_3(M) + \\
&+\frac{1}{4\pi} \int\limits_{\beta_{kn}}\, \Delta f_0(M)\left( \frac{1}{\rho_{M_0}(M)}-\frac{1}{\rho_{M_0}(M_{k,n-2})}\right)\, dm_3(M)+\\
&+\frac{1}{4\pi}\int\limits_{B_{C_1\Lambda_n}(M_{kn})\setminus\Omega_{n-2}^{\star}}\frac{\gamma_n \omega(\Lambda_n)}{\Lambda^{2}_n}\cdot \frac{1}{\rho_{M_0}(M_{k,n-2})}\, dm_3(M)+\\
&+\frac{1}{4\pi}\int\limits_{B_{C_1\Lambda_n}(M_{kn})\setminus\Omega_{n-2}^{\star}}\frac{\gamma_n \omega(\Lambda_n)}{\Lambda^{2}_n}\left( \frac{1}{\rho_{M_0}(M)}-\frac{1}{\rho_{M_0}(M_{k,n-2})}\right)\, dm_3(M)=\\
&=\frac{1}{4\pi}\cdot\frac{1}{\rho_{M_0}(M_{k,n-2})}\left(\, \int\limits_{\beta_{kn}}\, \Delta f_0(M)\, dm_3(M) + \int\limits_{\mathbb{R}^3}\, \phi_{kn}(M)\, dm_3(M) \right)+\\
&+\frac{1}{4\pi} \int\limits_{\beta_{kn}}\, \Delta f_0(M)\left( \frac{1}{\rho_{M_0}(M)}-\frac{1}{\rho_{M_0}(M_{k,n-2})}\right)\, dm_3(M)+\\
&+\frac{1}{4\pi}\int\limits_{B_{C_1\Lambda_n}(M_{kn})\setminus\Omega_{n-2}^{\star}}\frac{\gamma_n \omega(\Lambda_n)}{\Lambda^{2}_{n}}\left( \frac{1}{\rho_{M_0}(M)}-\frac{1}{\rho_{M_0}(M_{k,n-2})}\right)\, dm_3(M)=\\
&=\frac{1}{4\pi} \int\limits_{\beta_{kn}}\, \Delta f_0(M)\left( \frac{1}{\rho_{M_0}(M)}-\frac{1}{\rho_{M_0}(M_{k,n-2})}\right)\, dm_3(M)+\\
&+\frac{1}{4\pi}\int\limits_{B_{C_1\Lambda_n}(M_{kn})\setminus\Omega_{n-2}^{\star}}\, \frac{\gamma_n \omega(\Lambda_n)}{\Lambda^{2}_{n}}\left( \frac{1}{\rho_{M_0}(M)}-\frac{1}{\rho_{M_0}(M_{k,n-2})}\right)\, dm_3(M) = \\
&=A_k+D_k.
\end{split}
\end{equation}
We take into account that, for the indices $k$ in question and $M\in \beta_{kn}$, we have
\begin{equation}
\begin{split}
&\left|\frac{1}{\rho_{M_0}(M)}-\frac{1}{\rho_{M_0}(M_{k,n-2})}\right| =\left|\frac{1}{\|M_0 M\|}-\frac{1}{\|M_0 M_{k,n-2}\|}\right|\leq \\
&\leq\frac{C\Lambda_{n-2}}{{\|M_0 M_{k,n-2}\|}^2} \leq \frac{C\Lambda_{n-2}}{{|k-k_0|}^2 \Lambda^{2}_{n-2}}\leq\frac{C}{\Lambda_n{|k-k_0|}^2}.
\end{split}
\end{equation}
Since $d(M)\ge 2^{-n+1}$ for $M\in B_{C_1\Lambda_{n}}(M_{k,n-2})\setminus\Omega_{n-2}^{\star}$,   inequality (71) is also valid for such $M$ with a different $C$ depending on $C_1$ and $C_0$. Thus, due to (71) we get the following bounds for $A_k$ and $D_k$:

\begin{equation}
\begin{split}
&\left| A_k\right|\leq C\int\limits_{\beta_{kn}}\,\frac{\omega(d(M))}{d^{2}(M)}\cdot\frac{1}{\Lambda_n {|k-k_0|}^2}\, dm_3(M)\leq\\
&\leq C\omega(2^{-n})\cdot 2^{-n}\cdot \frac{1}{\Lambda_n {|k-k_0|}^2}\leq C\frac{\omega (2^{-n})}{(k-k_0)^2},
\end{split}
\end{equation}

\begin{equation}
\begin{split}
&\left| D_k\right|\leq C\int\limits_{B_{C_1\Lambda_n} (M_{k,n-2})\setminus\Omega_{n-2}^{\star}}\,\frac{\omega(\Lambda_{n})}{\Lambda^{2}_{n}}\cdot\frac{1}{\Lambda_n (k-k_0)^2}\, dm_3(M)\leq\\
&\leq C\frac{\omega (2^{-n})}{(k-k_0)^2}.
\end{split}
\end{equation}

\noindent
Consequently, (70), (72), and (73) imply
\begin{equation}
\begin{split}
&|\Sigma_1|+|\Sigma_3|\leq\sum_{k\leq k_0-2 \atop \text{ or } k\ge k_0+3}\, |A_k|+\sum_{k\leq k_0-2 \atop \text{ or } k\ge k_0+2} \, |D_k|\leq\\
&\leq C\omega(2^{-n})\sum_{\nu =1}^{\infty}\frac{1}{\nu^{2}}\leq C\omega(2^{-n}).
\end{split}
\end{equation}
Using (65)--(69) and (74), we have
\begin{equation}
\left| \upsilon_{2^{-n}}(M_0)-f(M_0)\right|\leq C\omega(2^{-n}).
\end{equation}
To get the required estimate (2) for any $\delta > 0$, we choose $n$ such that $2^{-n-1}<\delta \leq 2^{-n}$ and put $\upsilon_{\delta}=\upsilon_{2^{-n}}$;   relation (75) is equivalent to (2).

To verify estimate (3), we begin with the case $\delta =2^{-n}$. Let $\upsilon_{2^{-n}}$ be as before and let $M_0\in \Omega_{2^{-n}}(L)$. We have
\begin{equation*}
\begin{split}
&\left( \upsilon_{2^{-n}}(M_0)\right)'_{\bar{\nu}}=\frac{1}{4\pi}\int\limits_{\mathbb{R}^3\setminus\Omega_{n-2}^{\star}}\, \frac{\left(\rho_{M_1}(M)\right)'_{\bar{\nu}\mid_{M_1=M_0}}}{\rho^2_{M_0}(M)}\, \Delta f_0(M)\, dm_3(M) -\\
&-\frac{1}{4\pi}\int\limits_{\mathbb{R}^3\setminus\Omega_{n-2}^{\star}}\, \frac{\left(\rho_{M_1}(M)\right)'_{\bar{\nu}\mid_{M_1=M_0}}}{\rho^2_{M_0}(M)}\, \Phi_n(M)\, dm_3(M),
\end{split}
\end{equation*}
where
\noindent
$\bar{\nu}$ is an arbitrary unit vector. Then we get
\begin{equation}
\begin{split}
&\left|\left( \upsilon_{2^{-n}}(M_0)\right)'_{\bar{\nu}}\right|\leq C\int\limits_{\mathbb{R}^3\setminus\Omega_{n-2}^{\star}}\, \frac{\omega(d(M))}{\rho^2_{M_0}(M) d^2(M)}\, dm_3(M)+ \\
&+ C\int\limits_{\mathbb{R}^3\setminus\Omega_{n-2}^{\star}}\, \frac{|\Phi_{n}(M)|}{\rho^2_{M_0}(M)}\, dm_3(M) =\\
&= C\int\limits_{\left(\mathbb{R}^3\setminus\Omega_{n-2}^{\star}\right)\cap B_{2^{-n+3}|\Lambda|}(M_0)}\, \frac{\omega(d(M))}{\rho^2_{M_0}(M) d^2(M)}\, dm_3(M) + \\
&+ C\sum_{k=1}^{\infty}\int\limits_{\left( B_{2^{-n+k+3}|\Lambda|}(M_0)\setminus B_{2^{-n+k+2}|\Lambda|}(M_0)\right)\setminus\Omega_{n-2}^{\star}}\, \frac{\omega(d(M))}{\rho^2_{M_0}(M) d^2(M)}\, dm_3(M)+\\
&+ C\int\limits_{B_{2^{-n}\cdot |\Lambda^3|}(M_0)\setminus\Omega_{n-2}^{\star}}\, \frac{|\Phi_{n}(M)|}{\rho^2_{M_0}(M)}\, dm_3(M)+\\
&+ C\sum_{k=1}^{\infty}\int\limits_{\left( B_{2^{-n+k+3}|\Lambda|}(M_0)\setminus B_{2^{-n+k+2}|\Lambda|}(M_0)\right)\setminus\Omega_{n-2}^{\star}}\, \frac{|\Phi_{n}(M)|}{\rho^2_{M_0}(M)}\, dm_3(M).
\end{split}
\end{equation}
Due to (11)  we have  $d(M)\ge 2^{-n+1}$ for $M\notin \Omega_{n-2}^{\star}$. Hence  $$\rho_{M_0}(M)=\|M_0 M\|\ge 2^{-n+1}-2^{-n}=2^{-n}$$ for $M_0\in \Omega_{2^{-n}}(L)$. As in (46)--(52), we obtain
\begin{equation}
\begin{split}
&\int\limits_{B_{2^{-n+3}\cdot |\Lambda|}(M_0)\setminus\Omega_{n-2}^{\star}}\, \frac{\omega (d(M))}{\rho^{2}_{M_0}(M) d^2(M)}\, dm_3(M)\leq\\
&\leq \frac{C}{2^{-2n}}\int\limits_{B_{2^{-n+3}\cdot |\Lambda |}(M_0)}\, \frac{\omega (d(M))}{d^2(M)}\, dm_3(M)\leq\\
&\leq C\cdot 2^{2n}\cdot 2^{-n+3}\cdot |\Lambda |\cdot \omega (2^{-n+3}|\Lambda |)\leq C\cdot 2^{n}\omega (2^{-n}).
\end{split}
\end{equation}
Using the definition (61) of   $\phi_{kn}$ and the definition (63) of   $\Phi_n$, we obtain analogously that
\begin{equation}
\begin{split}
&\int\limits_{B_{2^{-n+3}\cdot |\Lambda |}(M_0)\setminus\Omega_{n-2}^{\star}}\, \frac{|\Phi_n(M)|}{\rho^2_{M_0}(M)}\, dm_3(M)\leq\\
&\leq C\cdot 2^{2n}\int\limits_{B_{2^{-n+3}\cdot |\Lambda |}(M_0)}\, |\Phi_n(M)|\, dm_3(M)\leq\\
&\leq C\cdot 2^{n}\omega (2^{-n}).
\end{split}
\end{equation}

Using analogs of   (46)--(52) once again, we get the estimates

\begin{equation}
\begin{split}
&\sum_{k=1}^{\infty}\int\limits_{\left( B_{2^{-n+k+3}|\Lambda |}(M_0)\setminus B_{2^{-n+k+2}|\Lambda |}(M_0)\right)\setminus\Omega_{n-2}^{\star}}\, \frac{\omega (d(M))}{\rho^2_{M_0}(M) d^2(M)}\, dm_3(M)\leq\\
&\leq C\sum_{k=1}^{\infty} 2^{2n-2k}\int\limits_{B_{2^{-n+k+3}\cdot |\Lambda |}(M_0)}\, \frac{\omega (d(M))}{d^2(M)}\, dm_3(M)\leq\\
&\leq C\sum_{k=1}^{\infty} 2^{2n-2k}\cdot 2^{-n+k}\cdot \omega (2^{-n+k})= C\cdot 2^{n}\sum_{k=1}^{\infty} 2^{-k}\omega (2^{-n+k})\leq\\
&\leq C\cdot 2^{n}\cdot \omega (2^{-n}).
\end{split}
\end{equation}
The last inequality in (79) is a consequence of the second part of assumption (1) concerning $\omega (t)$.
The definition (61) of   $\phi_{kn}$ allows us to deal with the function $|\Phi_n(M)|$ in the same way as with the expression $\frac{\omega (d(M))}{d^2(M)}$, so we get the relation
\begin{equation}
\sum_{k=1}^{\infty}\int\limits_{\left( B_{2^{-n+k+3}|\Lambda |}(M_0)\setminus B_{2^{-n+k+2}|\Lambda |}(M_0)\right)\setminus\Omega_{n-2}^{\star}}\, \frac{|\Phi_n(M)|}{\rho^2_{M_0}(M)}\, dm_3(M)\leq C\cdot 2^{n}\cdot \omega (2^{-n}).
\end{equation}
similar to (79).
Combining estimates (76)--(80), we come to the inequality
\begin{equation}
\left| \left( \upsilon_{2^{-n}}(M_0)\right)'_{\bar{\nu}}\right| \leq C\cdot 2^{n}\cdot \omega (2^{-n}).
\end{equation}
\noindent
This proves statement (3) for $\delta =2^{-n}$ since the constant $C$ in (81) is independent of $\bar{\nu}$. The case of  arbitrary $\delta >0$ is obtained in the same way as in the proof of statement~(2).

\subsection{Proof of Theorem~2}
\label{subsec:Theorem2}

We put
\begin{equation}
f^{\star}_0(x)=\int\limits_{0}^{x} \frac{\omega (t)}{t}\, dt,\quad x\in [0,\,1],
\end{equation}
and $f^{\star}_0(-x)=f^{\star}_0(x)$.
\noindent
Then   condition (1) implies that $f^{\star}_0(x)\leq C'\omega (x)$ and
\begin{equation}
f^{\star}_0(x)\ge \int\limits_{\frac{x}{2}}^{x} \frac{\omega (t)}{t}\, dt\ge \omega \left( \frac{x}{2}\right) \log{2} \ge\widetilde{C}'\omega (x), \quad x\in (0,\,1],
\end{equation}
where
\noindent
$\widetilde{C}'>0$ is independent of $x\in (0,\,1]$.
We have
\begin{equation}
f^{\star'}_0(x)=\frac{\omega (x)}{x}, \quad x\in (0,\,1].
\end{equation}

\noindent
Relations (1), (82)--(84) imply that $f^{\star}_0\in H^{\omega} ([0,\,1])$ and $f^{\star}_0(x)\asymp \omega (x)$.
We define $f_0(M)\stackrel{\rm def}= f^{\star}_0(x)$ for $M=(x,\,0,\,0)$. For $A>1$ and $0<x<\frac{1}{A}$, we have
\begin{equation*}
f^{\star}_0(A\, x)>\int\limits_{x}^{A\, x}\, \frac{\omega (t)}{t}\, dt \ge \omega (x) \log{A},
\end{equation*}
\noindent
and so

\begin{equation}
\omega (x)\leq \frac{1}{\log{A}} f^{\star}_0(A\, x) \leq \frac{\tilde{C}\omega (A\, x)}{\log{A}}
\end{equation}
Suppose there exist a sequence $\{k_\ell\}_{\ell =1}^{\infty}$ for which conditions (4) and (5) are fulfilled with some constants $C'_1$ and $C'_2$. We may assume that $\lambda_{k_\ell}>4$ for all $\ell$. Every function $V'_{k_{\ell} x}(M)$ is harmonic in the domain $\Omega_{\lambda_{k_\ell}\delta_{k_\ell}}([A_0,\, B_0])$, and (5) gives the following estimate:
\begin{equation}
\left| V'_{k_{\ell} x}(M)\right|\leq C'_2 \frac{\omega(\delta_{k_\ell})}{\delta_{k_\ell}}, \; M\in \Omega_{\lambda_{k_\ell}\delta_{k_\ell}}([A_0,\, B_0]).
\end{equation}
Let $r_{\ell} = \frac{1}{2}\lambda_{k_\ell}\delta_{k_\ell}$, and $A_{\ell}=\sqrt{\frac{1}{2}\lambda_{k_\ell}}$. We can use the Poisson integral   representation of the function $V'_{k_{\ell} x}$  harmonic in the ball $B_{2r_{\ell}}({O})$,
\begin{equation}
V'_{k_{\ell} x}(M)=\frac{1}{4\pi r_\ell}\int\limits_{\partial B_{r_\ell}({O})}\, V'_{k_{\ell} x}(P)\frac{r^2_{\ell}-\|{O}M\|^2}{\|M P\|^3}\, dm_{2}(P),
\end{equation}
\noindent
where $M\in B_{r_\ell}({O})$, and $dm_{2}(P)$ denotes the two-dimensional Lebesgue measure on the sphere $\partial B_{r_\ell}({O})$. If $M=(x,\,0,\,0)$, $|x|\leq A_{\ell}\delta_{k_\ell}$, then differentiating the integral (87) with respect to $x$ and taking into account   estimate (86), we obtain the inequality
\begin{equation}
\left| V''_{k_{\ell} xx}(M)\right|\leq C'_3\cdot \frac{1}{r_\ell}\max_{P\in \partial B_{r_\ell}({O})}\left| V'_{k_{\ell} x}(P)\right|\leq C'_4\frac{\omega (\delta_{k_\ell})}{r_\ell\delta_{k_\ell}}\leq 2C'_4 \frac{\omega (\delta_{k_\ell})}{\lambda{k_\ell})\delta^{2}_{k_\ell}}.
\end{equation}
Let $x_\ell =A_\ell\delta_{k_\ell}$ and $V^{\star}_k(x)=V_k((x,\,0,\,0))$.  Then (88) implies
\begin{equation}
\begin{split}
&\left| V^{\star}_{k_\ell}(x_\ell)+V^{\star}_{k_\ell}(-x_\ell)-2\,V^{\star}_{k_\ell}(0)\right|\leq \max_{|x|\leq x_\ell} \left| V''_{k_{\ell}}(x)\right|\cdot x^2_{\ell}\leq\\
&2C'_4 \frac{\omega (\delta_{k_\ell})}{\lambda{k_\ell})\delta^{2}_{k_\ell}}\cdot A^2_{\ell}\delta^2_{k_\ell}=C'_4\omega(\delta_{k_\ell}).
\end{split}
\end{equation}
From inequality (4) and the definition of $f_0$, it follows that
\begin{equation}
\begin{split}
&\left| (f^{\star}_0(x_\ell)-V^{\star}_{k_\ell}(x_\ell)+(f^{\star}_0(-x_\ell)-V^{\star}_{k_\ell}(-x_\ell))-2\, (f^{\star}_0(0)-V^{\star}_{k_\ell}(0)))\right|\leq\\
&\leq 4C'_1\,\omega (\delta_{k_\ell}).
\end{split}
\end{equation}
Estimates (89) and (90) put together imply  that
\begin{equation}
\left| f^{\star}_0(x_\ell)+f^{\star}_0(-x_\ell)-2\,f^{\star}_0(0)\right|\leq (C'_4+4C'_1)\omega (\delta_{k_\ell}).
\end{equation}

\noindent
On the other hand, $f^{\star}$ is an even function, so using relations (82) and (83), we get
\begin{equation}
\begin{split}
&f^{\star}_0(x_\ell)+f^{\star}_0(-x_\ell)-2\,f^{\star}_0(0)=2(f^{\star}_0(x_\ell)-f^{\star}_0(0))= \\
&= 2f^{\star}_0(x_\ell)\ge 2\widetilde{C'}\omega (x_\ell)=2\widetilde{C'}\omega (A_\ell\delta_{k_\ell}).
\end{split}
\end{equation}
\noindent
From (91) and (92), we obtain the inequality

\begin{equation}
2\widetilde{C'}\omega (A_\ell\delta_{k_\ell})\leq (C'_4+4C'_1)\omega (\delta_{k_\ell}).
\end{equation}

\noindent
Since
 $A_{\ell}\longrightarrow\infty$ as ${\ell\rightarrow \infty}$ and  inequality (93) is fulfilled for all $\ell$,   we have a contradiction with  inequality (85).  Theorem~2 is proved.

%



\section*{References}

\bibliography{mybib}

\end{document}